\documentclass[11pt]{amsart}
\usepackage{amsmath,amsthm, dsfont, amscd, amssymb, amsfonts,color}
\usepackage{hhline}
\usepackage[active]{srcltx}
\usepackage[all]{xy}
\usepackage{enumitem}

\newcommand{\com}{\textcolor{black}}

\newcommand{\comm}[1]{\textcolor{black}{\textbf{#1}}}

\newcommand\GK{\operatorname{GKdim}}

\newcommand\cB{\mathcal{B}}
\newcommand{\W}{{\mathcal W}}
\newcommand{\ku}{\Bbbk}
\newcommand{\N}{{\mathbb N}}
\newcommand{\U}{{\mathcal U}}
\newcommand{\ma}[1]{\hspace{-3pt}_{#1}\mathcal{M}}
\newcommand{\B}{{\mathcal B}}
\newcommand{\Vc}{{\mathcal V}}
\newcommand{\Ac}{{\mathcal A}}
\newcommand{\Ta}{{\mathcal T}}

\newcommand{\End}{\operatorname{End}}

\newcommand{\nuc}{\operatorname{ker}}
\newcommand{\la}{\langle}
\newcommand{\ra}{\rangle}
\newcommand{\gl}{{\mathfrak{gl}}}

\numberwithin{equation}{section}\theoremstyle{plain}

\newtheorem{theorem}{Theorem}[section]
\newtheorem{lema}[theorem]{Lemma}
\newtheorem{cor}[theorem]{Corollary}

\newtheorem{prop}[theorem]{Proposition}

\theoremstyle{definition}

\theoremstyle{remark}
\newtheorem{obs}[theorem]{Remark}

\newcommand\id{\operatorname{id}}

\newcommand\Ig{\operatorname{Im}}
\newcommand\ext{\operatorname{Ext}}

\def\pf{\begin{proof}}
\def\epf{\end{proof}}

\theoremstyle{remark}

\begin{document}

\thispagestyle{empty}

\title[Representations of the super Jordan plane]{Representations of the super Jordan plane}
\author[Andruskiewitsch, Bagio, Della Flora, Fl\^ores]
{Nicol\'as Andruskiewitsch, Dirceu Bagio, Saradia Della Flora, Daiana Fl\^ores}

\address{FaMAF-Universidad Nacional de C\'ordoba, CIEM (CONICET),
Medina A\-llen\-de s/n, Ciudad Universitaria, 5000 C\' ordoba, Rep\'ublica Argentina.} \email{andrus@famaf.unc.edu.ar}

\address{Departamento de Matem\'atica, Universidade Federal de Santa Maria,
97105-900, Santa Maria, RS, Brazil} \email{bagio@smail.ufsm.br, saradia.flora@ufsm.br, flores@ufsm.br}

\thanks{
N. A. was partially supported by CONICET, Secyt (UNC) and the MathAmSud project GR2HOPF. \\
D. B. and D. F. were supported by FAPERGS 2193-25.51/13-3, MATHAMSUD}

\dedicatory{To Antonio Paques.}

\begin{abstract}
It is shown that the finite-dimensional simple representations of the super Jordan plane $\cB$ are one-dimensional. The indecomposable representations of dimension $2$ and $3$ of $\cB$ are classified. Two families of indecomposable representations of $\cB$ of arbitrary dimension are presented. 
\end{abstract}

\maketitle

\section{Introduction}

Nichols algebras are graded connected algebras with a comultiplication in a braided sense.  In particular, the Jordan plane and the super Jordan plane are two Nichols algebras that
play an important role in the classification of pointed Hopf algebras with finite Gelfand-Kirillov dimension  \cite{aah-jordan, aah-triang}.

The Jordan plane was first defined in \cite{Gv} and considered in many papers, e.g. \cite{artin-schelter}, see also the references in \cite{aah-triang, Na}. Its representation theory was studied in \cite{Na}.

The purpose of this note is to begin the study of the representation theory of the super Jordan plane $\cB$: 
we classify the simple finite-dimensional $\cB$-modules (all of dimension 1, Theorem \ref{prop:classif-simples}) and the indecomposable 
$\cB$-modules of dimension 2 (Theorem \ref{theore:indec-2}) and 3 (Theorem \ref{theore:indec-3}).
We also observe that one of the generators of $\cB$ has at most two eigenvalues in every indecomposable $\cB$-module (Theorem \ref{teo:eigen}) and describe two families of indecomposable modules in every dimension.

\section{Basic facts}

\subsection{Notations and conventions} \label{subsec:notation}
Fix an algebraically closed field $\ku$ of characteristic 0; all vector spaces, tensor products, Hom spaces, algebras are over $\ku$.
All algebras are associative and all modules are left, unless explicitly stated. Let $A$ be a $\ku$-algebra; then $[ \, ,\, ]$ denotes the Lie bracket given by the commutator. As customary we use indistinctly the languages of modules and representations.
Denote by $\ma{A}$ the category of finite dimensional $A$-modules.
%and by $\Irr A$, respectively $\Inde A$, the set of isomorphism classes of irreducible, respectively indecomposable, $A$-modules. 
Given a $\ku$-vector space $V$, $\gl(V)$ denotes the Lie algebra of all linear operators on $V$. \com{The Jacobson radical of an algebra $A$ it will be denoted by $\operatorname{Jac}$ A.}

\subsection{The Jordan plane}\label{subsec:jordan}

The \emph{Jordan plane} is the free associative algebra $\mathcal{A}$  in generators $y_1$ and $y_2$ subject to the quadratic relation
\begin{align*}
 y_1y_2&-y_2y_1-y_2^{2}.
\end{align*}
The algebra $\mathcal{A}$ is a Nichols algebra,	$\GK \mathcal A = 2$ and $\{y_1^a y_2^b:  a,b\in\N_0\}$ is a basis of $\mathcal A$.  \com{By Proposition 3.4 of \cite{Na}, $\mathcal{A}$ is a Koszul algebra.}

\subsection{The super Jordan plane}\label{subsec:superjordan}

\smallbreak 
Let $x_{21}=x_1x_2+x_2x_1$ in the free associative algebra in generators $x_1$ and $x_2$. Let $\mathcal{B}$ be the algebra generated by $x_1$ and $x_2$ with defining relations
\begin{align}\label{eq:rels-B(V(-1,2))-1}
&x_1^2, \\
\label{eq:rels-B(V(-1,2))-2}
&x_2x_{21}- x_{21}x_2 - x_1x_{21}.
\end{align}
The algebra $\mathcal{B}$ (which is graded by $\deg x_1 = \deg x_2 =1$) was introduced in \cite{aah-jordan, aah-triang} and is called the \emph{super Jordan plane}. 
\com{Since $\B$ is not a quadratic algebra, it follows that $\B$ is not Koszul; see e.~g. $\S\,2.1$  of \cite{PP}. }

\begin{prop} \label{pr:sjp} \cite{aah-triang} The algebra $\cB$ is a Nichols algebra,	$\GK \cB = 2$ and $\{x_1^a x_{21}^bx_2^c: a\in\{0,1\},
	b,c\in\N_0\}$ is a basis of $\cB$. \qed
\end{prop}

The following identities are valid in $\B$:
\begin{align}\label{eq:relations1}
x_{21}x_1&=x_1x_{21},
\\\label{eq:relations2} x_2^2x_1 &=x_1x_2^2+x_1x_2x_1,
\\\label{eq:relations3} x_{21}x_2^2 &=(x_2^2-x_{21})x_{21}.
\end{align}
Indeed, in presence of \eqref{eq:rels-B(V(-1,2))-1}, \eqref{eq:rels-B(V(-1,2))-2} is
equivalent to \eqref{eq:relations2}.

By \eqref{eq:relations3} and Proposition \ref{pr:sjp}, the subalgebra of the super Jordan plane $\cB$ generated by $x_{2}^2$ and $x_{21}$, is isomorphic to the Jordan plane via $y_1 \mapsto x_2^{2}$ and $y_2 \mapsto x_{21}$.

\smallbreak

It is convenient to introduce $s = x_{21}$ and $t = x_2^2$.  By \eqref{eq:relations3}, $st=ts-s^2$ and whence
\begin{align}\label{eq:relations4}
[t, s^n] &= n s^{n+1},\,  n\geq 1;& x_1s &= sx_1; &
 x_2t &= tx_2; & tx_1 &= x_1(t + s).
\end{align}

\begin{lema}\label{lema:bc} Given $b,c \in \N$, we have that
	\begin{align*}
	& x_{21}^bx_2^c \overset{\ast}{=} (x_2-bx_1)x_{21}^bx_2^{c-1},& & x_1x_{21}^bx_2^c \overset{\heartsuit}{=}  x_1x_2x_{21}^bx_2^{c-1}.&
	\end{align*}
\end{lema}

\pf We prove $\ast$ by induction. For $b=c=1$, the relation is valid by \eqref{eq:rels-B(V(-1,2))-2}. Suppose that $\ast$ is valid for $b-1 > 0$ and $c=1$. Then
\begin{align*}
x_{21}^{b}x_2 &=x_{21}^{b-1}x_{21}x_2=x_{21}^{b-1}(x_2x_{21}-x_1x_{21})\\&=(x_2-(b-1)x_1)x_{21}^{b-1}x_{21}-x_{21}^{b}x_{1} =(x_2-bx_1)x_{21}^b.
\end{align*}
Fix $b\in \N$ and assume that the relation is true for $c-1$, with $c>1$. Thus
\[x_{21}^{b}x_2^c =(x_{21}^{b}x_2^{c-1})x_2=(x_2-bx_1)x_{21}^{b}x_2^{c-2}x_2=(x_2-bx_1)x_{21}^{b}x_2^{c-1}.\]
The proof of $\heartsuit$ is similar. \epf

The next result follows immediately from Proposition \ref{pr:sjp} 
and  Lemma \ref{lema:bc}.
\begin{prop}\label{cor:gener} The set $\{1,x_1,x_2,x_1x_2\}$ generates $\cB$ as a right $\Ac$-module. \qed
\end{prop}

%\com{Let $R\subset S$ be an extension of rings with the same unity. We recall that $S$ is a finite normalizing extension of $R$ if there exist $s_1,\ldots,s_n\in S$ such that $S=Rs_1+\ldots+ Rs_n$ and $Rs_i=s_iR$, for all $1\leq i\leq n$.}
%
%\com{
%\begin{cor}\label{cor:normal} $\Ac\subset \B$ is a finite normalizing extension.
%\end{cor}
%\pf By Proposition \ref{cor:gener}, $\B=\Ac 1+\Ac x_1+\Ac x_2+\Ac x_1x_2$.
%\epf}

\subsection{Simple modules}\label{subsec:reps}
Let $(V, \rho)$ be a finite-dimensional representation of $\B$; set $X_1 = \rho (x_{1})$, $X_2 = \rho (x_{2})$, $S= \rho (s)$ and $T = \rho (t)$ and 
\begin{align*}
V_0 = \ker X_1.
\end{align*}
 Then $V_0$ is always $\neq 0$ and it is stable under  $S$ and $T$ by \eqref{eq:relations4}.
In fact, let $E_{12}(n) \in \gl(\ku^n)$ (or $E_{12}$ if $n$ is clearly from de context) the matrix whose the entry $1\times 2$ is equal to $1$ and all other entries are equal to $0$.
Then the Jordan form of $X_1$ consists of $r$ blocks like $E_{12}(2)$ and $s$ blocks of size 1 filled by 0. Hence $\dim V = 2r+s$; $r=0 \iff V = V_0$.

\begin{lema}\label{kernel} Assume the previous notations. Then:
\begin{enumerate}[leftmargin=*,label=\rm{(\roman*)}]
\item\label{item:kernel1} $S$ and $T$ have a simultaneous eigenvector in $V_0$.
\item\label{item:kernel2}  $W =  X_{2}V_0 \cap V_0$ is a  submodule of $V$.
\item\label{item:kernel3} $U=X_{2} V_0 + V_0$ is a  submodule of $V$.
\end{enumerate}
\end{lema}

\pf \ref{item:kernel1}: The subspace of $\gl(V)$ generated by $T$ and $S^n$, $n \in \mathds{N}_0$, is a solvable Lie subalgebra by \eqref{eq:relations4}; then Lie Theorem applies.

\ref{item:kernel2}:  Clearly $X_1 W \subseteq X_1  V_0= \{0\} \subseteq W$. It remains to show that $X_2  W \subseteq W$. In fact, let $w \in W$, this is, $w \in V_0$ and  $w=X_2 v$ for some $v \in V_0$. Clearly $X_2 w \in X_2  V_0$. Moreover,
\begin{align*}
X_1  (X_2 w) & =  X_1X_2^2  v \overset{\eqref{eq:relations2}}{=} (X_2^2X_1-X_1X_2X_1 )  v = 0 \implies X_2  w \in W.
\end{align*}

\ref{item:kernel3}: Since $\B\cdot V_0\subseteq X_2 V_0$ and $\B\cdot (X_2 V_0)\subseteq V_0$,  the claim follows.
\epf

\begin{lema}\label{simples}
If $V \in \, \ma{\B}$ is simple,  then $V=V_0$.
\end{lema}
\pf
Assume that $V\neq V_0$. By Lemma \ref{kernel} we have that $W =  X_{2}  V_0 \cap V_0 = 0$ and $V =  X_{2}  V_0 + V_0$, so that $V =  X_{2}  V_0 \oplus V_0$.
By  Lemma \ref{kernel} \ref{item:kernel1}, there exists a simultaneous eigenvector $v \in V_0$ of $S$ and $T$, i.e. there exist $\alpha, \tau\in \ku$ such that $S v=\alpha v$, $T v=\tau v$.

Now $M=\text{span} \{v,X_2 v\} \neq 0$ is a $\B$-submodule of $V$ and by simplicity of $V$,
$M=V$. By our assumption, $X_2 v\notin V_0$; hence $\Lambda = \{v,X_2 v\}$ is a basis of $V$. Note that $[X_1]_{\Lambda}=\begin{pmatrix} 0 & \alpha \\
 0 & 0 \end{pmatrix}
$ and $[X_2]_{\Lambda} = \begin{pmatrix}
0 & \tau \\ 1 & 0 \end{pmatrix}
$. The relation $X_2^2X_1=X_1X_2^2+X_1X_2X_1$ is satisfied if and only if $\tau\alpha =\alpha\tau +\alpha^2$. Therefore $\alpha=0$ and  $V=V_0$, a contradiction.
\epf

Let  $A\in \End(\ku^n)$. Denote by $\ku^n_A$ the $\B$-module defined by   $X_1=0$ and $X_2=A$.
Every $\B$-module $V$ with $V = V_0$ is isomorphic to $\ku^n_A$ for some $A$. If  $B \in \End(\ku^m)$, then  $\ku^n_A\simeq \ku^m_B$ iff $n=m$ and $A$ and $B$ are similar matrices.

\begin{theorem}\label{prop:classif-simples} Every simple $\B$-module is isomorphic to $\ku^1_a$ for a unique $a\in \ku$.
\end{theorem}

\pf This follows from Lemma \ref{simples} and the preceding considerations.
\epf

\comm{
\begin{cor}\label{cor:jac} Let $\rho:\B\to \End V$ a finite dimensional representation of $\B$ and $B=\rho(\B)$. Then  there exists an integer $s$ such that $B/\!\operatorname{Jac} B\simeq \ku^s$ and $\operatorname{Jac} B=\{x\in B\,:\,x \text{ is nilpotent}\}$.
\end{cor}
\pf Since $B/\!\operatorname{Jac} B$ is semisimple and $\ku$ is algebraically closed, there are positive integers $n_1,\ldots,n_s$ such that $B/\!\operatorname{Jac} B=M_{n_1}(\ku)\times\cdots \times M_{n_s}(\ku)$. The composition  
\[\xymatrix{& \B\ar@{->>}[r]^{\rho}  & B \ar@{->>}[r]^-\pi& B/\!\operatorname{Jac} B  \ar@{->>}[r]^-{\pi_j}&  M_{n_j}(\ku) }\]
is a finite dimensional simple representation of $\B$. Hence, by Theorem \eqref{prop:classif-simples}, $n_1=\cdots=n_s=1$. Thus, $B/\!\operatorname{Jac} B\simeq \ku^s$.  
Let $x\in B$ a nilpotent element. Then $\pi(x)$ is a nilpotent element of $B/\!\operatorname{Jac} B$. Since $B/\!\operatorname{Jac} B$ is commutative, we obtain that $\pi(x)\in \operatorname{Jac} \left(B/\!\operatorname{Jac} B\right)=\{0\}$. Hence, $x\in \operatorname{Jac} B$. On the other hand, $B$ finite dimensional implies that $\operatorname{Jac} B$ is a nilpotent ideal. Consequently, $\operatorname{Jac} B=\{x\in B\,:\,x \text{ is nilpotent}\}$.
\epf}

%\simeq 

We also remark:

\begin{prop}\label{prop:indec-x1trivial} If $V$ is an indecomposable $\B$-module with $V= V_0$, then there exist $n\in \N$ and $\lambda \in \ku$ such that $V$ is isomorphic to $\ku^n_A$ where $A$ is the Jordan block of size $n$ with  eigenvalue $\lambda$. \qed
\end{prop}

If $A$ is the Jordan block of size $n$ with  eigenvalue $\lambda$, then denote $\Ac_\lambda = \ku^n_A$.

\subsection{Indecomposable modules}
Throughout this subsection, $V$, $X_1$, $X_2$, $T$ and $S$ are as in $\S\,$\ref{subsec:reps}.
When $V$ is indecomposable, we will prove that $T$ has a unique eigenvalue. In order to do this, the following relations are useful.
\begin{lema}\label{lem:other-rel} Let $\lambda \in \ku$, $z:=t-\lambda \id  \in \cB$ and $n\in \N$. Then
	\begin{align*}
	&z^nx_1 \overset{\clubsuit}{=}
	x_1\sum\limits_{j=0}^{n}\frac{n!}{(n-j)!}s^jz^{n-j},& & z^{n}x_1x_2
\overset{\diamondsuit}{=}	x_1x_2\sum\limits_{j=0}^{n}\frac{n!}{(n-j)!}s^jz^{n-j}. &
	\end{align*}	
	
\end{lema}
\pf We prove $\clubsuit$ by induction on $n$;  the proof of $\diamondsuit$ is similar. We will use that  $x_1zs^n=x_1x_{21}^nz+nx_1s^{n+1}$, which can be verified easily. Note that \[zx_1=x_1z+x_1x_2x_1=x_1z+x_1s=x_1(z+s),\]
and whence the formula is true for $n=1$. Denote $\zeta_{n,j}:=\frac{n!}{(n-j)!}$, $0\leq j\leq n$. Consider $n>1$ and assume that the formula is true for $n-1$. Then
\begin{align*}
z^nx_1=(zx_1)\sum\limits_{j=0}^{n-1}\zeta_{n-1,j}s^jz^{n-1-j}    \overset{\eqref{eq:rels-B(V(-1,2))-2}}{=}(x_1z+x_1s)\sum\limits_{j=0}^{n-1}\zeta_{n-1,j}s^jz^{n-1-j}\\
      =\sum\limits_{j=0}^{n-1}\zeta_{n-1,j}(x_1zs^j)z^{n-1-j}+\sum\limits_{j=0}^{n-1}\zeta_{n-1,j}x_1s^{j+1}z^{n-1-j}
      =x_1\sum\limits_{j=0}^{n}\zeta_{n,j}s^jz^{n-j}.
\end{align*}
\epf

Let $\lambda$ be an eigenvalue of $T$. Denote by $V_{\lambda}^T$ the generalized eigenspace of $V$ associated to $\lambda$, i.~e.  $V_{\lambda}^T:=\cup_{j\geq 0} \nuc\left(T-\lambda\id\right)^j$

\begin{lema}\label{lema:componente} $V_{\lambda}^T$ is a $\cB$-submodule of $V$, for all eigenvalue $\lambda$ of $T$.
\end{lema}
\pf  Clearly $V_{\lambda}^{T}=\nuc\left(T-\lambda\id\right)^r=\nuc\left(X_2^2-\lambda\id\right)^r$, where $r$ is the maximal size of $\lambda$-blocks in the Jordan normal form of $T$. Thus $V_{\lambda}^{T}$ is stable by $X_2$. It remains to show that it stable by  $X_1$. By Lemma \ref{lem:other-rel}, if $u\in V_{\lambda}^{T}$ then
\begin{align*}
(T-\lambda\id)^nX_1 u= X_1\sum\limits_{j=0}^{n}\zeta_{j,n}S^j(T-\lambda\id)^{n-j}u.
\end{align*}
By Lemma 2.1 of \cite{Na}, $S$ is nilpotent. Taking $n$ big enough, it follows that $(T-\lambda\id)^nX_1 u=0$ and whence $X_1u\in V_{\lambda}^T$. 
\epf

Now  Lemma \ref{lema:componente} implies the next result.

\begin{theorem}\label{teo:eigen} Let $\lambda_1,\ldots,\lambda_t$ be the different eigenvalues of $T$. Then $V$ decomposes into the direct sum of the $\cB$-submodules $V_{\lambda_i}^T$.

In particular, if $V$ is indecomposable then $T$ has a unique eigenvalue. Hence either $X_2$ has a unique eigenvalue or  else
the eigenvalues of $X_2$ are $\lambda$ and $-\lambda$, with $\lambda\in \ku^{\times}$. \qed
\end{theorem}

\com{ Given $\lambda\in \ku$, denote by $\,\,\ma{\B}_{\lambda}$ the full subcategory of $\,\,\ma{\B}$ whose objects are the $\B$-modules $V$ such that $V=\ker\, (T-\lambda\id)^m$, for some $m\in \N_0$. With this notation, the next result follows immediately from Theorem \ref {teo:eigen}.
\begin{cor} $\ma{\B}\simeq \prod_{\lambda\in \ku}\,\,\ma{\B}_{\lambda}$. \qed	
\end{cor}}

The next result will be useful in $\S\,$3.

\begin{lema}\label{dim1}
Let $\Lambda =\{v_1,\cdots, v_n\}$ be a basis of $V$ such that $[X_1]_{\Lambda} =E_{12}$ and $W$ a one-dimensional $\B$-submodule of $V$. Then:
	\begin{enumerate}[leftmargin=*,label=\rm{(\roman*)}]
\item\label{item:dim11} If $L$ is a complement  (as a $\B$-module) of $W$ in $V$  then  $ L\cap V_0 =\langle v_1\rangle$.
\item\label{item:dim12}  $W= \langle v_1 \rangle$ does not have a complement (as a $\B$-submodule) in $V$. \end{enumerate}
\end{lema}
\pf \ref{item:dim11}: Assume that $W=\langle w\rangle$ and $\{u_1,u_2, \cdots, u_{n-1}\}$ is a basis of $L$. Since $W_0=W\cap V_0 \neq 0$, it follows that $W\subset V_0$.  Using that $v_2$ is a linear combination of $w,u_1,u_2, \cdots, u_{n-1}$ we see that $v_1=X_1 v_2\in L$; hence $v_1\in V_0\cap L$. \smallbreak

\ref{item:dim12}: It follows at once from \ref{item:dim11}.
\epf

\section{Indecomposable representations of dimension $2$ and $3$}\label{subsec:indec2}

\subsection{Dimension $2$}
In this subsection we describe all $2$-dimensional indecomposable representations of $\B$. Fix $(V,\rho)$ a $2$-dimensional representation of $\cB$.

\begin{lema}
If $V \neq V_0$ then $V$ is indecomposable.
\end{lema}
\pf Suppose that $V$ is decomposable, i.e. there are non-trivial submodules $U$ and $W$  such that $V= U \oplus W$. Then
$V_0= U_0 \oplus W_0 =  U \oplus W = V$.
\epf

Define representations of $\B$  on the vector space $\ku^2$ given by $X_1 = E_{12}$ and the following action of $x_2$:
\begin{itemize}[leftmargin=*]\renewcommand{\labelitemi}{$\diamond$}
	\item   $X_2=\begin{pmatrix}
	a & b \\
	0 &  a \end{pmatrix}$,  $a,b\in \ku$.  This is denoted by $\U_{a,b}$.
	
	\item  $X_2=\begin{pmatrix}
	a & 0 \\
	0 & - a \end{pmatrix}$, $a\in \ku^{\times}$. This is denoted by  $\Vc_a$.
	
\end{itemize}

It is easy to check that these are indecomposable modules pairwise non-isomorphic.

\begin{theorem}\label{theore:indec-2}
	Every 2-dimensional indecomposable representation  of $\B$ is isomorphic either to $\U_{a,b}$, or to $\Vc_a$, or to $\ku^2_{\lambda}$ for unique $a,b, \lambda\in \ku$.
\end{theorem}

This confirms Theorem \ref{teo:eigen}.

\pf If $V = V_0$, then Proposition \ref{prop:indec-x1trivial} applies. Assume that $V_0 \neq 0$; then there exists a basis $\Lambda =\{v_1,v_2\}$ of $V$ such that $[X_1]_{\Lambda}= E_{12}$. Let $[X_2]_{\Lambda} =\begin{pmatrix} a & b \\   c & d              \end{pmatrix}$.
Then \eqref{eq:relations2} is satisfied if and only if
\begin{align*} c(a+d)&=0& &\text{and} & d^2+c &= a^2.
\end{align*}

\noindent Suppose that $c\neq 0$. Then by the first equation it follows that $d=-a$. Replacing in the second equation we have that $c=0$, which is a contradiction. Therefore $c=0$ and consequently $d=a$ or $d=-a$.

If $d = a$ then $V\simeq \U_{a,b}$. Assume that $d = -a \neq 0$ and take $w_1=v_1$ and $w_2=\frac{-b}{2a} v_1+v_2$. Then $\Omega =\{w_1,w_2\}$ is a basis of $V$ such that $[X_1]_{\Omega}= E_{12}$ and  $[X_2]_{\Omega} =\begin{pmatrix}
a & 0 \\
0 & - a \end{pmatrix}$. Thus $V\simeq \Vc_a$.
\epf

\begin{cor}\label{cor:ext}
	If $\ext^1(\ku_a^1,\ku_b^1)\neq 0$ then $a=\pm b$. \qed
\end{cor}

\subsection{Dimension $3$}\label{subsec:indec3}
Let $V$ be a $\B$-module of dimension 3 such that $V\neq V_0$. Throughout this subsection, $\Lambda =\{v_1,v_2, v_3\}$ denotes a basis of $V$ such that $[X_1]_{\Lambda} =E_{12}$. We define four families of representations of $\B$ on the vector space $V$ determined by the following action of $[X_2]_{\Lambda}$, for all $a,b,c,d,e\in \ku$:
\begin{align*}
&\Theta_1 :  \begin{pmatrix}a & b  & c\\
\noalign{\medskip}       0 & d  & e \\
\noalign{\medskip}       0 & \frac{a^2-d^2}{e}  & -d\\
\end{pmatrix},\,\,e\in \ku^{\times};&
&\Theta_2 :  \begin{pmatrix}a & b  & c\\
\noalign{\medskip}       0 & a  & 0 \\
\noalign{\medskip}       0 & d  & e\\
\end{pmatrix};&\\
& \Theta_3 : \begin{pmatrix}  a & b  &\frac{c^2-a^2}{d}\\
\noalign{\medskip}        0 & c  & 0\\
\noalign{\medskip}        d & e  & -a\\
\end{pmatrix},\,\, d\in \ku^{\times};&
& \Theta_4 : \begin{pmatrix}  a & b  & c\\
\noalign{\medskip}        0 & -a & 0\\
\noalign{\medskip}        0 & d  & e\\
\end{pmatrix},\,\, a \in \ku^{\times}.&
\end{align*}

\begin{lema}
The families $\Theta_1,\, \Theta_2,\, \Theta_3$ and $\Theta_4$ contain all $3$-dimensional representations of $\B$, up to isomorphism.
\end{lema}
\pf Let $[X_2]_{\Lambda} =
\begin{pmatrix}
\alpha & \beta & \gamma  \\
\delta & \epsilon & \zeta \\
\eta & \theta & \iota
\end{pmatrix}$. Then \eqref{eq:relations2} is valid if and only if
\begin{align}\label{sistema}
\left\{\begin{array}{l}
\delta(\alpha+\epsilon)=-\zeta\eta \\
\zeta(\epsilon+\iota)=-\gamma\delta \\
\eta(\alpha+\iota)=-\delta\theta\\
\epsilon^2-\alpha^2= -\delta +\gamma\eta-\zeta\theta
\end{array}\right.
\end{align}

\noindent {\it Claim:} If the system \eqref{sistema} has solution then $\delta=0$. \smallbreak

\noindent Assume that $\delta \neq 0$. If $\zeta=0$ then $\gamma=0$ and $\epsilon=-\alpha$. Thus, $\delta=0$ which is a contradiction. If $\zeta\neq 0$ then \[\gamma=\frac{-\zeta(\epsilon+\iota)}{\delta}, \quad \eta=\frac{-\delta(\alpha+\epsilon)}{\zeta}\quad \text{and}\quad \theta=\frac{(\alpha+\epsilon)(\alpha+\iota)}{\zeta}.\]
From the last equation of \eqref{sistema}, $\delta=0$ which is again a contradiction. \smallbreak

\noindent Assume $\delta=0$. Thus $\zeta\eta=0$. If $\zeta\neq 0$ then $\eta=0$, $\iota=-\epsilon$ and $\theta=\frac{\alpha^2-\epsilon^2}{\zeta}$. Hence $V$ belongs to the family $\Theta_1$. When $\zeta=0$ and $\eta\neq 0$, it follows that $\iota=-\alpha$ and $\gamma=\frac{\epsilon^2-\alpha^2}{\eta}$. Thus, $V$ belongs to  the family $\Theta_3$. If $\zeta=0$ and $\eta=0$ then $\epsilon=|\alpha|$. In this case, $V$ belongs to the families $\Theta_2$ or $\Theta_4$.\epf

\begin{obs}\label{generator} Let $L$ a $\B$-submodule of $V$ of dimension $2$ such that $L\cap V_0$ is one-dimensional. Fix $\overline{L}:=L/(L\cap V_0)=\langle \overline{u} \rangle$. Since $u\notin V_0$, we can suppose that $u=\alpha v_1+ v_2+\gamma v_3\in L$, with $\alpha,\,\gamma\in \ku$.
\end{obs}

\begin{prop}\label{prop:indecom} Let $V$ be a $\B$-module. Then:
\begin{enumerate}[leftmargin=*,label=\rm{(\roman*)}]
\item\label{item:indecom1} the representations in the family $\Theta_1$ are always indecomposable;	
\item\label{item:indecom2} a representation  in the family $\Theta_2$ is indecomposable if and only if $c \neq 0$ and $e=a$ or  $d \neq 0$ and $e=a$;	
\item\label{item:indecom3} the representations in the family  $\Theta_3$ are always indecomposable;	
\item\label{item:indecom4} a representation  in the family  $\Theta_4$ is indecomposable if and only if $c \neq 0$ and $e=a$ or $d \neq 0$ and $e=-a$.
\end{enumerate}	
\end{prop}

\pf\ref{item:indecom1}: The unique one-dimensional $\B$-submodule of $V$ is $\la v_1 \ra$ which does not have complement by Lemma \ref{dim1} (ii). \smallbreak

\ref{item:indecom2}: Let $V$ be a representation of $\B$ of the type $\Theta_2$. Suppose that $W= \langle w \rangle$  is a one-dimensional $\B$-submodule of $V$. Since $W\subset V_0$, see $\S\,\ref{subsec:reps}$, $w=\alpha v_1 +\beta v_3$, with $\alpha,\beta\in \ku$. Note that $X_2 w=\gamma w$, $\gamma\in \ku$, if and only if $\beta(\gamma-e)=0$ and $\alpha(\gamma-a)=\beta c$. Consequently, the one-dimensional $\B$-submodules of $V$ are:
\begin{itemize}[leftmargin=*]\renewcommand{\labelitemi}{$\diamond$}
\item $\la v_1 \ra$, $\la v_3 \ra$, $c=0$, $e\neq a$,
\item $\la \alpha v_1+\beta v_3 \ra$, $c=0$, $e=a$,
\item $\la v_1 \ra$, $\la v_1+\frac{e-a}{c}v_3 \ra$, $c\neq 0$, $e\neq a$
\item $\la v_1\ra$, $c \neq 0$, $e=a$.
\end{itemize}

\noindent Assume $e \neq a$. If $c \neq 0$, $V=\la v_1+\frac{e-a}{c}v_3 \ra\oplus \la v_1, v_1+v_2+\frac{d}{a-e}v_3 \ra$. If $c=0$,  $V=\la v_3\ra \oplus \langle v_1,v_1+v_2+\frac{d}{a-e}v_3 \rangle$. If $c=d=0$, $V=\la v_1,v_2\ra\oplus \la v_3\ra$. Hence, $V$ is decomposable.

\noindent Conversely, suppose  $e=a$ and $c \neq 0$. Then the unique one-dimensional $\B$-submodule of $V$ is $\la v_1\ra$ which does not have complement. Suppose that $e=a$ and $d \neq 0$.
Assume that $W$ is a one-dimensional $\B$-submodule of $V$ which admits a complement $L=\la u_1, u_2 \ra$. Then by Lemma \ref{dim1} (ii), $v_1 \in L \cap V_0$. By Remark \ref{generator}, $\overline{L}= \la \overline{u} \ra$ where $u= \alpha v_1+v_2+\beta v_3$. Thus $X_2 \overline{u}= \gamma \overline{u}$, $\gamma \in \ku$, if and only if $\gamma=a$ and $\beta (a-e)=d$. Since $d\neq 0$ and $e=a$ then $W$ does not have complement in $V$.

\ref{item:indecom3}: Suppose that $W= \langle w \rangle$  is a one-dimensional $\B$-submodule of $V$ which admits a complement $L$. Then by Lemma \ref{dim1} (ii), $ \la v_1 \ra = L \cap V_0$, which is a contradiction because $d \neq 0$.

\ref{item:indecom4}: Analogous to item (ii).
\epf

\subsubsection{Isomorphism classes in $\Theta_1$}\label{susubsec:noniso-theta1}
Assume $V$ in the family $\Theta_1$. We distinguish: for all  $a,b,c,d,e\in \ku$

\begin{itemize}[leftmargin=*]\renewcommand{\labelitemi}{$\diamond$}
	\item $X_2=\begin{pmatrix}
	a & b & c \\
	0 &  d & e\\ 0 & \frac{a^2-d^2}{e} & -d \end{pmatrix}, \,\, e \in \ku^{\times}$. This is denoted by $\mathcal{Y}_{a,b,c,d,e}.$\smallbreak
	\item  $X_2=\begin{pmatrix}
	a & b & 0 \\
	0 &  a & 1\\ 0 & 0 & -a \end{pmatrix}$. This is denoted by $\mathcal{U}^{a,b}.$
\end{itemize}
\noindent By Proposition \ref{prop:indecom} (i), these representations are indecomposable. Note that $\mathcal{U}^{a,b}=\mathcal{Y}_{a,b,0,a,1}$.

% Nicolás: como \mathcal{U}^{a,b}=\mathcal{Y}_{a,b,0,a,1}, no hace falta mirar los otros casos de isomorfismos en la proxima proposition.

\begin{prop}\label{prop-non-iso-family3}
Every $3$-dimensional indecomposable representation $V$ of $\B$ in $\Theta_1$ is isomorphic either to $\mathcal{U}^{a,b}$, or to $\mathcal{Y}_{a,b,c,d,e}$. Moreover,
\begin{align*}
\mathcal{Y}_{a,b,c,d,e} \simeq \mathcal{Y}_{a,b',c',d',e'}\text{ if and only if } (a-d')\displaystyle\frac{ce'-c'e}{e'}=e(b'-b)+c(d'-d).
\end{align*}
	 
In particular, $\mathcal{U}^{a,b}\simeq \mathcal{U}^{a,b'}$ if and only if $b=b'$.
\end{prop}		
\pf  Since $\la X_2v_1\ra=\Ig X_1$, we obtain that $a$ is invariant. Consider the indecomposable representation $\mathcal{Y}_{a,b',c',d',e'}$ of $\B$.
If $d'=a$, taking the basis $\{v_1,\frac{c'}{e'}v_1+v_2, \frac{1}{e'}v_3\}$ we conclude that $\mathcal{Y}_{a,b',c',d',e'} \simeq \mathcal{U}^{a,b'}$.\smallbreak

\noindent Note that $\mathcal{Y}_{a,b,c,d,e}$ and $\mathcal{Y}_{a,b',c,d',e'}$ are isomorphic if and only if there exists a basis $\{w_1,w_2,w_3\}$ of $V$ such that $X_1w_1=X_1 w_3=0$, $X_1w_2=w_1$,  $X_2 w_1=aw_1$,  $X_2 w_2= b'w_1+d'w_2+ \frac{a^2-{d'}^{2}}{e'}w_3$ and  $X_2  w_3= c'w_1+e'w_2-d'w_3$. Since $\la v_1\ra=\Ig X_1$  and $V_0$ has dimension $2$, then we can consider $w_1=v_1$, $w_2=\lambda_1v_1+\lambda_2v_2+\lambda_3v_3$ and $w_3=\beta_1v_1+\beta_3v_3$, $\,\lambda_1,\lambda_2,\lambda_3,\beta_1,\beta_3 \in \ku$. Then, $\mathcal{Y}_{a,b,c,d,e}\simeq \mathcal{Y}_{a,b',c,d',e'}$  if and only if $(a-d')\displaystyle\frac{ce'-c'e}{e'}=e(b'-b)+c(d'-d)$.  
\epf

\subsubsection{Isomorphism classes in $\Theta_2$}\label{susubsec:noniso-theta2} Consider $V$ in the family $\Theta_2$ and
the following distinguish representations: for all $a\in \ku$
\begin{itemize}[leftmargin=*]\renewcommand{\labelitemi}{$\diamond$}
	\item   $X_2=\begin{pmatrix}
	a & 0 & 0 \\
	0 &  a & 0\\ 0 & 1 & a \end{pmatrix}$.  This is denoted by $\mathcal{R}_{a}$. \smallbreak
	
	\item  $X_2=\begin{pmatrix}
	a & 0 & 1 \\
	0 &  a & 0 \\ 0 & 0 & a\end{pmatrix}$. This is denoted by $\mathcal{S}_{a}$.\smallbreak
	
		\item  $X_2=\begin{pmatrix}
		a & 0 & b \\
		0 &  a & 0 \\ 0 & c & a\end{pmatrix},\,\, b\in \ku^{\times}$ or $c \in \ku^{\times}$. This is denoted by $\mathcal{T}_{a,b,c}$.
\end{itemize}

\noindent By Proposition \ref{prop:indecom} (ii),  these are indecomposable representations. Notice that $\mathcal{R}_{a}=\mathcal{T}_{a,0,1}$ and $\mathcal{S}_{a}=\mathcal{T}_{a,1,0}$.

\begin{prop}\label{prop-non-iso-family2}
Every $3$-dimensional indecomposable representation $V$ of $\B$  in $\Theta_2$ is isomorphic either to $\mathcal{R}_{a}$, or to $\mathcal{S}_{a}$ or to $\mathcal{T}_{a,b,c}$. Moreover,
$\Ta_{a,b,c}$ and $\Ta_{a,b',c'}$ are isomorphic if and only if $bc=b'c'$.
\end{prop}		
\pf Let  $V'$ the representation of $\B$ given by
\begin{align*}
 [X_2]_{\Lambda}= \begin{pmatrix}
a & d' & b' \\
0 &  a & 0 \\ 0 & c' & a\end{pmatrix}.
\end{align*}
If $b'=0$, then  by Proposition \ref{prop:indecom} (ii) we have that $c' \neq 0$. In this case, taking the basis $\{v_1,v_2,d'v_1+c'v_3\}$ of $V'$, we conclude that $V' \simeq \mathcal{R}_{a}$.
Similarly, if $c'=0$ then $b' \neq 0$. Taking the basis $\{v_1,v_2-\frac{d'}{b'}v_3,\frac{1}{b'}v_3\}$ of $V'$, we obtain that $V' \simeq \mathcal{S}_{a}$. If $b,b',c,c'\in \ku^{\times}$, taking the basis $\{v_1,v_2,\frac{d'}{c'}v_1+v_3\}$, it follows that $V'\simeq \Ta_{a,b',c'}$.

Finally, notice that $\Ta_{a,b,c}\simeq \Ta_{a,b',c'}$ if and only if there exists a basis $\{w_1,w_2,w_3\}$ of $\ku^3$ such that $X_1 w_1=X_1 w_3=0$, $X_1 w_2=w_1$,   $X_2w_1=aw_1$,  $X_2w_2= aw_2+cw_3$ and  $X_2w_3= bw_1+aw_3$. We can assume $w_1=v_1$, $w_2=\lambda_1v_1+\lambda_2v_2+\lambda_3v_3$ and $w_3=\beta_1v_1+\beta_3v_3$, $\,\lambda_1,\lambda_2,\lambda_3,\beta_1,\beta_3 \in \ku$ . Note that $X_2 w_2=w_1$ if and only if $\lambda_2=1$. Moreover, $X_2w_2=aw_2+cw_3$ and $X_2 w_3=bw_1+aw_3$ if and only if $bc=b'c'$. \epf

\subsubsection{Isomorphism classes in $\Theta_3$}\label{susubsec:noniso-theta3}
Consider $V$ in the family $\Theta_3$ and  the following distinguished representations: for all $a,b,c,d,e\in \ku $

 \begin{itemize}[leftmargin=*]\renewcommand{\labelitemi}{$\diamond$}
  \item $X_2=\begin{pmatrix}
	a & b & \frac{c^2-a^2}{d} \\
	0 &  c & 0\\ d & e & -a \end{pmatrix},\,\, d\in \ku^{\times}$. This is denoted by $\mathcal{W}_{a,b,c,d,e}$. \smallbreak
	\item  $X_2=\begin{pmatrix}
	a & b & 0 \\
	0 &  a & 0\\ 1 & 0 & -a \end{pmatrix}$. This is denoted by $\mathcal{U}_{a,b}$.
\end{itemize}
\noindent By Proposition \ref{prop:indecom} (iii), these representations are indecomposable. Observe that $\mathcal{U}_{a,b}=\mathcal{W}_{a,b,a,1,0}$.

\begin{prop}\label{prop-non-iso-family3bis}
Every $3$-dimensional indecomposable representation $V$ of $\B$  in $\Theta_3$ is isomorphic either to $\mathcal{U}_{a,b}$, or to $\mathcal{W}_{a,b,c,d,e}$. Moreover, 
\begin{align*}
\mathcal{W}_{a,b,c,d,e}\simeq \mathcal{W}_{a',b',c,d',e'} \text{ if and only if } \frac{ae-bd-ce}{d}=\frac{a'e'-b'd'-ce'}{d'}.
\end{align*}
In particular, $\mathcal{U}_{a,b}\simeq \mathcal{U}_{a,b'}$ iff $b=b'$.
\end{prop}		
\pf  Since the characteristic polynomial of $X_2$ is $(t-c)^2(t+c)$, $c$ is an invariant. Let the indecomposable representation $\mathcal{W}_{a',b',c,d',e'}$ of $\cB$. If $c=a$, taking the basis $\{v_1,-\frac{e}{d}v_1+v_2, dv_3\}$ we conclude that $\mathcal{W}_{a',b',c,d',e'} \simeq \mathcal{U}_{a',b'}$. Note that $\mathcal{W}_{a,b,c,d,e}\simeq \mathcal{W}_{a',b',c,d',e'}$ if and only if there is a basis $\{w_1,w_2,w_3\}$ of $\ku^3$ such that $X_1 w_1=X_1 w_3=0$, $X_1 w_2=w_1$,  $X_2 w_1=aw_1$,  $x_2w_2= aw_2+cw_3$ and  $X_2 w_3= bw_1+aw_3$. We can assume $w_1=v_1$, $w_2=\lambda_1v_1+\lambda_2v_2+\lambda_3v_3$ and $w_3=\beta_1v_1+\beta_3v_3$, where $\lambda_1,\lambda_2,\lambda_3,\beta_1,\beta_3 \in \ku$. However $X_2 w_1=a'w_1+d'w_3$ if and only if $\beta_1=\frac{a-a'}{d'}$ e $\beta_3=\frac{d}{d'}$. With this choose of $\beta_1$ and $\beta_3$ we have that $X_2 w_3=\frac{c^2-a'^{2}}{d'}w_1-a'w_3$. Finally, $X_2 w_2=b_1w_1+cw_2+e'w_3$ if and only if $\frac{ae-bd-ce}{d}=\frac{a'e'-b'd'-ce'}{d'}$. \epf

\subsubsection{Isomorphism classes in $\Theta_4$}\label{susubsec:noniso-theta4}
Consider $V$ in $\Theta_4$ and the following distinguish representations: for all $a\in \ku^{\times}$
\begin{itemize}[leftmargin=*]\renewcommand{\labelitemi}{$\diamond$}
	\item   $X_2=\begin{pmatrix}
	a & 0 & 1 \\
	0 &  -a & 0\\ 0 & 0 & a \end{pmatrix}$.  This is denoted by $\mathcal{V}^{a}$.\smallbreak
	
	\item  $X_2=\begin{pmatrix}
	a & 0 & 0\\
	0 &  -a & 0 \\ 0 & 1 &- a\end{pmatrix}$. This is denoted by $\mathcal{V}_{a}$.
\end{itemize}

\noindent By Proposition \ref{prop:indecom} (iv), these are indecomposable representations pairwise non-isomorphic.

\begin{prop}\label{prop-non-iso-family4}
	Every $3$-dimensional indecomposable representation $V$ of $\B$  in $\Theta_4$ is isomorphic either to $\mathcal{V}^{a}$ or to $\mathcal{V}_{a}$ for unique $a \in \ku^{\times}$.
\end{prop}		
\pf Let $V'$ be a $3$-dimensional indecomposable representation of $\B$ such that
\begin{align*}
[X_2]_{\Lambda}= \begin{pmatrix}
a & b & c \\
0 &  -a & 0 \\ 0 & d & e\end{pmatrix}, \,\,a\in \ku^{\times}.
\end{align*}
Since $V'$ is indecomposable, by Proposition \ref{prop:indecom} (iv) we have that $c \neq 0$ and $e=a$ or $d \neq 0$ and $e=-a$. If $c \neq 0$ and $e=a$, taking the basis $\{v_1,\frac{cd-2ab}{4a^2}v_1+v_2-\frac{d}{2a}v_3,v_1+\frac{1}{c}v_3\}$ of $V'$, we obtain $V' \simeq \mathcal{V}^{a}$.
If $d\neq 0$ and $e = -a$, taking the basis $\{v_1,-\frac{2ab+cd}{4a^2}v_1+v_2, -\frac{dc}{2a}v_1+dv_3\}$ of $V'$, it follows that $V' \simeq \mathcal{V}_{a}$.
\epf

\subsubsection{Classification of indecomposable 3-dimensional $\cB$-modules}\label{susubsec:noniso-general}

\begin{theorem}\label{theore:indec-3}
	Every $3$-dimensional indecomposable $\cB$-module is isomorphic either to  $\ku^3_{\lambda}$ for a unique $\lambda$, or else to a representation in one of the families $\Theta_j$, $j=1,2,3,4$, with the constraints described in Proposition \ref{prop:indecom}. The isomorphism classes are described in Propositions \ref{prop-non-iso-family3}, \ref{prop-non-iso-family2}, \ref{prop-non-iso-family3bis} and \ref{prop-non-iso-family4}. \qed
\end{theorem}

Again, this agrees with Theorem \ref{teo:eigen}.

\begin{obs}
It is straightforward to verify that two $3$-dimensional indecomposable representations of $\cB$ that belong to different families $\Theta_i$, $i=1,2,3,4$, are not isomorphic.
\end{obs}

\section{Families of indecomposable $\B$-modules}\label{subsec:family}

Throughout this section $(V,\rho)$ is an $n$-dimensional representation of $\cB$, $\Lambda=\{v_1,\ldots,v_n\}$ is a basis of $V$, $X_1=\rho(x_1)$, $X_2=\rho(x_2)$ and $[X_1]_{\Lambda}=E_{12}$.

\subsection{The family \textbf{$\U_a$}}\label{subsec:4.1}
Let $a\in\ku$. Consider the following action of $X_2$ on $V$:
\begin{align*}
[X_2]_{\Lambda} =
\begin{pmatrix}
a & 0 & 0 & 0 &\ldots & 0& 0\\
0 & a & 0 & 0 & \ldots & 0& 0\\
0 & 1 & a & 0 & \ldots & 0 & 0\\
0 & 0 & 1 & a & \ldots & 0 & 0\\
\vdots &\vdots & \vdots &\ddots &\ddots  & \vdots& \vdots\\
0 & 0 & 0 & \ldots & 1 & a & 0\\
0 & 0 & 0 & \ldots &0 & 1 & a &
\end{pmatrix}.
\end{align*}
Clearly $V$ with this action is a $\B$-module which will be denoted by $\U_a$.

\begin{lema}\label{lem:zero}
Let $\W$ be a proper $\B$-submodule of $\U_a$. Then:
\begin{enumerate}[leftmargin=*,label=\rm{(\roman*)}]
	\item\label{item:zero1} $v_2\notin \W$;
	\item\label{item:zero2}  If  $v=\sum_{i=1}^n\lambda_iv_i\in \W$ then $\lambda_2=0$.
\end{enumerate}
\end{lema}
\pf \ref{item:zero1}: Suppose $v_2\in \W$. Then $v_1=X_1 v_2\in \W$ and $X_2 v_2=av_2+v_3\in \W$. Hence $v_3\in \W$. Again, $X_2 v_3=av_3+v_4\in \W$ and consequently $v_4\in \W$. With this procedure, we obtain that $\Lambda\subset \W$. Thus, $\W=\U_a$ and we have a contradiction. \smallbreak

\ref{item:zero2}: Assume $\lambda_2\neq 0$ and fix $w_1=\lambda_2^{-1}v$. Thus $w_1=\alpha_1 v_1+v_2+\ldots +\alpha_nv_n$, where $\alpha_i=\lambda_2^{-1}\lambda_i$, for all $1\leq i\leq n$. Consider the following elements of $V$:
\begin{align*}
 &w_j:=v_{j+1}+\alpha_3v_{j+2}+\ldots+\alpha_{n-j+1}v_n,\quad \text {for all }\, 2\leq j\leq n-2.&
\end{align*}
By a straightforward calculation, we obtain that $X_2 w_j=aw_j+w_{j+1}$, for all $1\leq j\leq n-2$. Thus, $w_1,\ldots,w_{n-2}\in \W$ and $X_2 w_{n-2}=aw_{n-2}+v_n$. Therefore, $v_n\in \W$. But  $w_{n-2}=v_{n-1}+\alpha_3v_n$ and whence $v_{n-1}\in \W$. By this procedure, it follows that $v_3,\ldots,v_n\in \W$. From $v_1=X_1 w_1 \in \W$, it follows that $v_2\in \W$ which contradicts {\rm (i)}. \epf

\begin{theorem}\label{teo:fam1}
	$\U_a$ is an indecomposable $\B$-module, for all  $n\geq 2$.
\end{theorem}
\pf Suppose $\U_a$ decomposable. Let $\W,\,\widetilde{\W}$ be nontrivial $\B$-submodules of $\U_a$ such that $\U_a=\W\oplus\widetilde{\W}$. Consider $\{w_1,\ldots,w_r\}$ and $\{w_{r+1},\ldots,w_n\}$ basis of $\W$ and $\widetilde{\W}$ respectively. By Lemma \ref{lem:zero}, $w_i=\lambda_{i1}v_1+\lambda_{i3}v_3+\ldots+\lambda_{in}v_n$, for all $1\leq i\leq n$. Since $v_2\in \U_{a}$, there exist $\alpha_1,\ldots,\alpha_n\in \ku$ such that $v_2=\alpha_1w_1+\ldots+\alpha_nw_n$, a contradiction.\epf

%\begin{obs}\label{obs:fam1}
%We denote by $A$ and $B$ respectively, the matrices $[X_1]_{\Lambda}$ and $[X_2]_{\Lambda}$ that define the $\B$-module $\U_a(n)$. Also, $\mathfrak{h}(A,B)$ denotes the subalgebra  of $\gl(V)$ generates by $A$ and $B$. Since $AB=BA=aA$, it follows that $\mathfrak{h}(A,B)=\ku\cdot A \oplus \ku\la B\ra$ as vector spaces and $\ku\cdot A$ is a two-sided ideal of $\mathfrak{h}(A,B)$.
%\end{obs}

\subsection{The family \textbf{$\Vc_a$}}\label{subsec:4.2}
Let $a\in\ku^{\times}$.  Consider the following action of $X_2$ on $V$:
\begin{align*}
[X_2]_{\Lambda} =
\begin{pmatrix}
a & 0 & 0 & 0 &\ldots & 0& 0\\
0 & -a & 0 & 0 & \ldots & 0& 0\\
0 & 1 & -a & 0 & \ldots &0&0\\
0 & 0 & 1 & -a & \ldots & 0& 0\\
\vdots &\vdots & \vdots &\ddots &\ddots  & \vdots& \vdots\\
0 & 0 & 0 & \ldots & 1 & -a & 0\\
0 & 0 & 0 & \ldots &0 & 1 & -a &
\end{pmatrix}.
\end{align*}

% Nicolás: acá no se puede cambiar el -a por un elemento b. En efecto, si uno hace la cuenta X_2^2X_1=X_1X_2^2+X_1X_2X_1 (reemplazando -a por b en X_2) sale que a^2=b^2, o sea, b=a o -a.

Notice that $V$ is a $\B$-module which will be denoted by $\Vc_a$. Since $a\neq 0$, $\U_a$ and $\Vc_a$ are not isomorphic.
\begin{theorem}\label{teo:fam2}
	$\Vc_a$ is an indecomposable $\B$-module, for all  $n\geq 2$.
\end{theorem}
\pf Let $\W$ a proper $\B$-submodule of $\Vc_a$. As in Lemma \ref{lem:zero} (i), we can show that $v_2\notin \W$. Let  $v\in \W$ such that $v=\sum_{i=1}^n\lambda_iv_i$. Assume that $\lambda_2\neq 0$ and consider $u:=\lambda_2^{-1}v\in \W$. Then $v_1=X_1 u\in \W$. Take $w_1:=u-\lambda_2^{-1}\lambda_1 v_1$ and note that $w_1= \alpha_2v_2+\ldots +\alpha_nv_n$, where $\alpha_i=\lambda_2^{-1}\lambda_i$, for all $2\leq i\leq n$. Considering the following elements of $V$
\begin{align*}
&w_j:=v_{j+1}+\alpha_3v_{j+2}+\ldots+\alpha_{n-j+1}v_n,\quad \text {for all }\, 2\leq j\leq n-2,&
\end{align*}
it follows $X_2 w_j=-aw_j+w_{j+1}$, for all $1\leq j\leq n-2$. As in Lemma \ref{lem:zero}, this implies that $v_2\in \W$ which is a contradiction. Thus, the result follows as in Theorem \ref{teo:fam1}. \epf

%\begin{obs}\label{obs:fam2}
%Denote by $A$ and $C$ respectively, the matrices $[X_1]_{\Lambda}$ and $[X_2]_{\Lambda}$ that define the $\B$-module $\Vc_a(n)$. As in Lemma \ref{obs:fam1}, $\mathfrak{h}(A,C)$ denotes the subalgebra  of $\gl(V)$ generates by $A$ and $C$. Since $CA=-AC=aA$, we obtain $\mathfrak{h}(A,C)=\ku\cdot A \oplus \ku\la C\ra$ and $\ku\cdot A$ is a two-sided ideal of $\mathfrak{h}(A,C)$.
%\end{obs}

\end{document}